\documentclass[5p,times]{elsarticle_preprint}
\usepackage{amsmath}
\usepackage{bm}
\usepackage{amssymb}
\usepackage{subfig}
\usepackage{enumitem}
\usepackage{comment}
\usepackage{graphicx}

\newtheorem{thm}{Theorem}
\newtheorem{prop}{Proposition}
\newtheorem{cor}{Corollary}

\newdefinition{rmk}{Remark}
\newproof{pf}{Proof}
\newproof{pot}{Proof of Theorem \ref{thm2}}

\begin{document}
	
\journal{Automatica}

\title{A recursively feasible and convergent Sequential Convex Programming procedure to solve non-convex  problems with linear equality constraints}

\author{Josep Virgili-Llop and Marcello Romano}
\address{Mechanical and Aerospace Engineering Department, Naval Postgraduate School, 1 University Circle, Moneterey, CA, 93940}

\begin{abstract}
A computationally efficient method to solve non-convex programming problems with linear equality constraints is presented. 
The proposed method is based on a recursively feasible and descending sequential convex programming procedure proven to converge to a locally optimal solution.
Assuming that the first convex problem in the sequence is feasible, these properties are obtained by convexifying the non-convex cost and inequality constraints with inner-convex approximations.
Additionally, a computationally efficient method is introduced to obtain inner-convex approximations based on Taylor series expansions.
These Taylor-based inner-convex approximations provide the overall algorithm with a quadratic rate of convergence.
The proposed method is capable of solving problems of practical interest in real-time.
This is illustrated with a numerical simulation of an aerial vehicle trajectory optimization problem on commercial-of-the-shelf embedded computers.
\end{abstract}


\maketitle

\section{Introduction}
The ability to solve trajectory optimization problems in real-time on embedded platforms is of paramount importance to achieve a high level of autonomy in aerospace applications. 
To accomplish this goal, obtaining assurances about an algorithm's computational feasibility, tractability, and strict convergence is a necessity \cite{Tsiotras2017}.
The growing interest in the computational aspects of trajectory optimization has recently crystallized with the coining of the term Computational Guidance and Control \cite{Lu2017b}.  

Convex optimization paired with solvers based on interior-point methods enable to solve convex programming problems to global optimality in polynomial-time \cite{Nesterov1994,Boyd2004}. These properties make convex programming suitable for implementation in real-time autonomous applications \cite{Dueri2016,Scharf2016,Virgili2018_IJRR}. 

Many optimization problems can be cast as convex programming problems \cite{Boyd2004}. For example, fixed-final-time trajectory optimization problems of linear systems with convex cost and inequality constraints form convex programming problems when discretized \cite{Hull1997}.
Some problems with apparently non-convex control constraints can also be cast as convex programming problems thanks to lossless convexification or exact relaxation techniques \cite{Accikmecse2011,Blackmore2012,Accikmecse2013,Harris2014,Liu2016,Verscheure2009}. In trajectory optimization, the use of convex programming techniques has been remarkably fruitful in the past decade \cite{Liu2017}. Some of the areas where it has been applied include: powered soft-landing \cite{Acikmese2007,Blackmore2010,Accikmecse2013,Sagliano2017,Pinson2018}, time-optimal path following \cite{Verscheure2009,Debrouwere2013,Lipp2014}, spacecraft proximity maneuvering \cite{Mueller2008,Harris2014b,Lu2013,Liu2014}, re-entry guidance \cite{Liu2015,Liu2015b,Zhao2017}, swarm control \cite{Alonso2015,Morgan2016}, quadrotor guidance \cite{Szmuk2017}, missile guidance \cite{Liu2016}, and motion planning \cite{Schulman2014,Misra2017,Virgili2018_IJRR}.

Despite these significant advances, there is a large number of optimization problems of important practical interest that cannot be cast as convex programming problems. For example, problems involving nonlinear systems or problems with non-convex constraints generally produce non-convex programming problems.
Although a wide range of approaches are available to solve these non-convex programming problems \cite{Nocedal2006}, none of these provide the same computationally attractive properties that interior-point methods offer for convex problems.
A notable exception may be the Lassere's relaxations \cite{Lasserre2001} and their variations \cite{Nie2013}. Lassere-type relaxations are able to relax non-convex programming problems, written in terms of polynomials, to semidefinite programming problems. Although promising, these relaxations dramatically increase the dimensionality of the problem.

In an attempt to bridge some of the attractive properties of convex programming to non-convex problems, the sequential convex programming (SCP) procedure  has been proposed \cite{Liu2017}. First developed to solve structural problems \cite{Svanberg1987,Zillober2001}, this method has been adopted by the trajectory optimization community (in addition to previous references see \cite{Dinh2010,Szmuk2016,Mao2017,Alrifaee2017,Szmuk2018,Mao2018}). In an SCP procedure, a convex approximation of the original non-convex problem is repeatedly solved until convergence. Finding the solution of a non-convex programming problem is thus reduced to solving a collection of convex---hence tractable---programming problems.

To obtain the sequence of convex programming problems that constitute an SCP, the non-convex functions are usually linearized around the previous iteration's solution, in what is known as successive linearization \cite{Liu2017,Mao2018}. Linearization-based SCP procedures can suffer from inadmissible iterates, infeasible problems, or cost increases between iterations. Trust regions and penalties can be used to alleviate these problems, and when properly applied can result in a convergent SCP procedure with a superlinear rate of convergence \cite{Mao2017,Mao2018}.%

\subsection{Contributions}
Here we advance the SCP's state-of-the-art by making two critical contributions. First, we present an SCP procedure that offers guaranteed convergence to a locally optimal solution for non-convex problems with linear equality constraints. This procedure is recursively feasible as long as the first convex problem in the sequence is feasible. Trust regions are not needed as the solution to the convexified problem is always admissible.
These advantageous properties are obtained when inner-convex approximations are used \cite{Marks1978}.
These properties were first demonstrated by Marks and Wright \cite{Marks1978} but they have remained largely unnoticed, and in some cases, independently re-discovered by the aerospace trajectory optimization community \cite{Liu2014,Morgan2016}. Here, we generalize and consolidate previously available results and highlight how the use of inner-convex approximations prevents inadmissible solutions, infeasible problems, and cost increases. Additionally, we briefly review the different available methods to obtain inner-convex approximations (namely, sum-of-squares decompositions and regularizations).

Secondly, we introduce a computationally efficient approach to obtain inner-convex approximations. Our approach is based on the convexification of the Taylor series expansion of a non-convex analytic function. 
These inner-convex approximations allow the SCP to achieve a locally quadric rate of convergence under mild conditions. As Taylor series may be infinite and therefore truncated, we introduce a high-order regularization approach that allows to truncate the Taylor series and retain the inner-convex properties.

The combination of a recursively feasible and convergent SCP with an efficient method to obtain inner-convex approximations results in an algorithm that is suitable for onboard implementation and real-time use. As each iteration involves solving a convex programming problem, the time per iteration is bounded. Then, the recursively feasible aspect allows the SCP to stop at any iteration and still yield an admissible solution. Finally, the convergent property ensures that each solution is better than the last. This properties make the proposed method particularly attractive for online trajectory optimization on autonomous systems.

\subsection{Organization of the paper}
The rest of this paper is organized as follows. First, the class of non-convex problems with linear equality constraints is introduced in Section~\ref{sec:Problem}. Section~\ref{sec:SCP} outlines the generic SCP procedure. A brief discussion on the successive linearization approach is included in Section~\ref{sec:Suc_linear}. In Section~\ref{sec:Descending} we introduce the concept of inner-convex approximations and present the recursively feasible and convergent SCP procedure that is obtained when these approximations are used.
In the following Section~\ref{sec:Methods} a brief overview of the existing methods that can be used to find inner-convex approximations is provided.
Then, in Section~\ref{sec:IC_Taylor} we introduce our method to obtain inner-convex approximations based on the convexification of a Taylor series expansion. Finally, and to illustrate the use of the proposed method, a numerical simulation example, based on the path-planning of an aerial vehicle, is provided in Section \ref{sec:Simulation}. The results of an implementation on two different embedded computers are also provided, substantiating the claim that the proposed method is suitable for onboard implementation and real-time use.

\section{Non-convex programming problems with linear equality constrains}
\label{sec:Problem}
The focus of this work is directed towards non-convex programming problems with linear equality constraints. In canonical form, this class of problems can be written as:

\textbf{Problem $P$. Non-convex programming problem with linear equality constraints}
\begin{align*}
\text{minimize: } & f_{0}\left(x\right) \\
\text{subject to: } & f_{i}\left(x\right)\leq 0, && \text{for }i=1\ldots p \tag{$P$}\label{prob:NonConvex}\\
& h_{i}\left(x\right)= 0, &&\text{for }i=1\ldots q
\end{align*}
In this problem, $x\in\mathbb{R}^{n}$ is the optimization variable, the equality constraint functions $h_{i}:\mathbb{R}^{n}\rightarrow\mathbb{R}$ are affine, and the cost or objective function $f_{0}:\mathbb{R}^{n}\rightarrow\mathbb{R}$ as well as the inequality constraint functions $f_{i}:\mathbb{R}^{n}\rightarrow\mathbb{R}$ are non-convex. To simplify notation, generic non-convex functions are simply denoted by $f:\mathbb{R}^{n}\rightarrow\mathbb{R}$ in the following.

In general, we call an admissible point of the programming problem $P$ any point that satisfies the constraints, $x\in\mathcal{D}_{P}$, with $\mathcal{D}_{P}$ denoting the problem's admissible set. The programming problem is said to be feasible if $\mathcal{D}_{P}$ has a non-empty interior and infeasible if $\mathcal{D}_{P}=\varnothing$.

A locally optimal solution of problem $P$, denoted by $x^{\star}_{P}$, is understood as a point where the Karush-Kuhn-Tucker (KKT) conditions are satisfied \cite{Boyd2004}:
\begin{subequations}
	\begin{align}
	f_{i}\left(x^{\star}_{P}\right)&\leq 0, && \text{for }i=1\ldots p \label{eq:primal_ineq}\\
	h_{i}\left(x^{\star}_{P}\right)&= 0, && \text{for }i=1\ldots q \label{eq:primal_eq}\\
	\nu_{P,i}^{\star}&\geq0,&& \text{for }i=1\ldots p \label{eq:dual}\\
	\nu_{P,i}^{\star}f_{i}\left(x^{\star}_{P}\right)&=0,&& \text{for }i=1\ldots p \label{eq:compl}
	\end{align}
	\begin{equation}
	\nabla f_{0}\left(x^{\star}_{P}\right)+\sum_{i=1}^{p}{\nu_{P,i}^{\star}\nabla f_{i}\left(x^{\star}_{P}\right)}+\sum_{i=1}^{q}{\mu_{P,i}^{\star}\nabla h_{i}\left(x^{\star}_{P}\right)}=0\label{eq:stationarity}
	\end{equation}
\end{subequations}
with $\nu_{P}^{\star}\in\mathbb{R}^{p}$ and $\mu_{P}^{\star}\in\mathbb{R}^{q}$ denoting the optimal dual variables associated with $x^{\star}_{P}$. Finally, $\nabla f$ is used to denote the gradient of a function $f$.

\section{Sequential Convex Programming (SCP) procedure}
\label{sec:SCP}
In a sequential convex programming procedure a locally optimal solution to the original non-convex problem $P$ is found by iteratively solving a convex approximation of the original non-convex problem. The iterations within the SCP's sequence are represented by $k\in\mathbb{N}$ and the solutions to each of the convex programs in the sequence are denoted by $x_{\left(k\right)}$. The $k$-th convex programming problem $P_{\mathrm{cvx}}^{\left(k\right)}$ is obtained by convexifying the non-convex functions, $f$, around the previous iteration's solution $x_{\left(k-1\right)}$. 

The notation $f^{\left(x_{e}\right)}_{\text{cvx}}$ is used here to denote a convex approximation of $f$ around the expansion point $x_{e}$. For simplicity, the approximation around the previous iteration's solution $f_{\mathrm{cvx}}^{\left(x_{\left(k-1\right)}\right)}$, is simply denoted as $f_{\mathrm{cvx}}^{\left(k-1\right)}$. With this notation, the convex problem to be repeatedly solved in an SCP procedure is:

\textbf{Problem $P_{\mathrm{cvx}}^{\left(k\right)}$. Convex approximation of problem $P$ around the previous iteration's solution $x_{\left(k-1\right)}$}
\begin{align*}
\text{minimize: } & f_{\text{cvx},0}^{\left(k-1\right)}\left(x\right) \\
\text{subject to: } & f_{\text{cvx},i}^{\left(k-1\right)}\left(x\right)\leq 0, && \text{for }i=1\ldots p \tag{$P_{\mathrm{cvx}}^{\left(k\right)}$}\label{prob:ConvexApprox}\\
& h_{i}\left(x\right)= 0, && \text{for }i=1\ldots q
\end{align*}
The first time problem $P_{\mathrm{cvx}}^{\left(k\right)}$ is solved, $k=1$, the problem is convexified around an initial guess, $x_{\left(0\right)}$. The SCP procedure iteratively solves problem $P_{\mathrm{cvx}}^{\left(k\right)}$ until a stopping criterion is satisfied. For example, when the cost difference between iterates falls below a certain threshold $\epsilon\ge 0$, \emph{i.e.},
\begin{equation}
\left\lvert f_{0}\left(x_{\left(k-1\right)}\right)-f_{0}\left(x_{\left(k\right)}\right)\right\rvert\leq\epsilon.
\end{equation}

\section{SCP with successive linearization}
\label{sec:Suc_linear}
A common approach to obtain the convex approximations $f_{\text{cvx}}^{\left(x_{e}\right)}$ is to linearize $f$ \cite{Liu2017},
\begin{equation}
f^{\left(x_{e}\right)}_{\text{cvx}}\left(x\right)=f\left(x_{e}\right)+\delta x^{T}\nabla f\left(x_{e}\right),
\end{equation}
\begin{equation}
\delta x=x-x_{e}.
\end{equation}
It has been shown that a linearization-based SCP is locally convergent \cite{Dinh2010}. However, the convergence is not guaranteed outside a small neighborhood around the locally optimal solution.
In general, a linearization-based SCP procedure may yield inadmissible solutions. \emph{I.e}, the solutions, $x_{\left(k\right)}$, when re-evaluated with the original non-convex functions, $f$, may violate the inequality constraints. Furthermore, the cost may be higher than the cost from previous iterations, \emph{i.e.}, $f_{0}\left(x_{\left(k\right)}\right)>f_{0}\left(x_{\left(k-1\right)}\right)$ \cite{Mao2017b}. This means that it is not known \emph{a priori} if solving the next convex problem in the sequence will yield a better or even an admissible solution. Additionally, linearizing the cost can lead to unbounded problems, completely preventing the use of a linearization-based SCP.

To mitigate these negative effects, trust regions are commonly used \cite{Conn2000}. Trust regions restrict the solution, $x_{\left(k\right)}$, to a closed bounded convex set, $x\in\mathcal{T}\subset\mathbb{R}^{n}$, where the linearization is expected to be accurate.
The heuristic use of trust regions does not  eliminate the possibility of an inadmissible solution or cost increase. If that solution can be discarded and the same convex problem re-solved with tighter trust regions, the SCP is shown to be convergent \cite{Mao2017}.

 
\subsection{Penalty method to resolve infeasibility}
When a convex problem in the sequence is infeasible a penalty method can be applied in an attempt to gain feasibility \cite{Schulman2014,Mao2017b}. The penalty method consists on turning the infeasible  inequality constraints into parts of the cost function. First, the inequality $f_{i}\left(x\right)\leq 0$ is relaxed with a slack variable $s_{i}$,
\begin{subequations}
	\begin{align}
	f_{i}\left(x\right)\leq s_{i}, \quad& \text{for }i=1\ldots l\leq p \\
	s_{i}\geq 0, \quad& \text{for }i=1\ldots l\leq p
	\end{align}
\end{subequations}
and then, the cost function is augmented as
\begin{equation}
f_{0}\left(x\right)+\kappa_{\left(k\right)}\sum_{i=1}^{l}s_{i},
\end{equation}
with $\kappa_{\left(k\right)}$ denoting the penalty's weight or gain.

The theory of exact penalties \cite{Pillo1994} states that the constrained and the penalty-relaxed problems are equivalent as long as $\kappa_{\left(k\right)}$ is greater than the largest optimal dual variable, $\nu_{P,i}^{\star}$, associated with the relaxed inequalities. The values of $\nu_{P,i}^{\star}$ are unlikely to be known \emph{a priori}, which prompts the use of large gains.

The use of penalties is a useful strategy to resolve infeasibilities, and a critical tool to find an admissible solution when an admissible initial guess is not available. The use of penalties also allows the SCP to proceed. However, the SCP may converge to an inadmissible point, $x_{\left(\infty\right)}\notin\mathcal{D}_{P}$.

\section{ SCP with inner-convex approximations}
\label{sec:Descending}
If the convex programing problems within the SCP are constructed with inner-convex approximations the SCP is recursively feasible, descending, and convergent \cite{Marks1978}. No trust regions are required and the solution is always improving as the cost decreases between iterations.

Inner-convex approximations are those convex approximations, $f_{\text{cvx}}^{\left(x_{e}\right)}$, that meet the three following conditions:
\begin{subequations}
	\label{eq:ConvexProp}
	\begin{equation}
	f_{\text{cvx}}^{\left(x_{e}\right)}\left(x\right)\geq f\left(x\right),
	\label{eq:ConvexProp_geq}
	\end{equation}
	\begin{equation}
	f_{\text{cvx}}^{\left(x_{e}\right)}\left(x_{e}\right)=f\left(x_{e}\right),
	\label{eq:ConvexProp_eq}
	\end{equation}
	\begin{equation}
	\nabla f_{\text{cvx}}^{\left(x_{e}\right)}\left(x_{e}\right)= \nabla f\left(x_{e}\right),
	\label{eq:ConvexProp_eqDerivative}
	\end{equation}
	and additionally, as a convex function, $f_{\text{cvx}}^{\left(x_{e}\right)}$, must also obey for $\tau\in\left[0,1\right]$ \cite{Boyd2004},
	\begin{equation}
	f_{\text{cvx}}^{\left(x_{e}\right)}\left(\tau x_{1}{+}\left(1{-}\tau\right)x_{2}\right){\geq} \tau f_{\text{cvx}}^{\left(x_{e}\right)}\left(x_{1}\right){+}\left(1{-}\tau\right)f_{\text{cvx}}^{\left(x_{e}\right)}\left(x_{2}\right)
	\label{eq:ConvexProp_cvx}.
	\end{equation}
\end{subequations}
For completeness we now introduce the theorem and proof of the SCP's recursive feasibility, descent, and convergence, first credited to Marks and Wright \cite{Marks1978}.
\begin{thm}
\label{th:main_theorem}
An SCP that uses inner-convex approximations, as defined in Eq.~\eqref{eq:ConvexProp}, is an infinite descent algorithm, $f_{0}\left(x_{\left(k\right)}\right)\leq f_{0}\left(x_{\left(k-1\right)}\right)$, guaranteeing convergence to a locally optimal solution, $x_{\left(\infty\right)}=x^{\star}_{P}$, and with all its iterates remaining admissible points of the original non-convex problem, $x_{\left(k\right)}\in\mathcal{D}_{P}$, as long as the first convex problem in the sequence $P_{\mathrm{cvx}}^{\left(1\right)}$ is feasible, $\mathcal{D}_{P_{\mathrm{cvx}}^{\left(1\right)}}\neq\varnothing$, and as long as the original non-convex problem~$P$ is lower bounded $f_{0}\left(x\right)>-\infty$.
\end{thm}
The proof of Theorem~\ref{th:main_theorem} is supported by the following two propositions:
\begin{prop}
	The admissible set of the convexified problem $P_{\mathrm{cvx}}^{\left(k\right)}$ is a subset of the original problem's $P$ admissible set $\mathcal{D}_{P_{\mathrm{cvx}}^{\left(k\right)}}\subseteq\mathcal{D}_{P}$.
	\label{prop:admissible_prop}
\end{prop}

\begin{pf}
	Invoking Eq.~\eqref{eq:ConvexProp_geq}
	\begin{equation}
	f_{i}\left(x\right) \leq f_{\text{cvx},i}^{\left(k-1\right)}\left(x\right),\quad \text{for }i=1\ldots p,
	\end{equation}
	shows that if a point $x$ meets the convexified constraints, $x\in\mathcal{D}_{P_{\mathrm{cvx}}^{\left(k\right)}}$, \emph{i.e.}, $f_{\text{cvx},i}^{\left(k-1\right)}\left(x\right)\leq 0$, it also meets the original non-convex constraints, \emph{i.e.}, $f_{i}\left(x\right)\leq 0$, hence $x\in\mathcal{D}_{P}$\qed
\end{pf}

\begin{cor}
	From Proposition~\ref{prop:admissible_prop} it follows that the solution to the convexified problem~$P_{\mathrm{cvx}}^{\left(k\right)}$ is an admissible point of the original non-convex problem~$P$: $x_{\left(k\right)}\in\mathcal{D}_{P}$.
	\label{prop:feasible}
\end{cor}
\begin{prop}
	The convexified problem~$P_{\mathrm{cvx}}^{\left(k\right)}$ is feasible if an admissible point to the original non-convex problem~$P$ is used as the expansion point of the convex approximation, $f_{\text{cvx},i}^{\left(k-1\right)}$.
	\label{prop:admissible}
\end{prop}
\begin{pf}
	Given Eq.~\eqref{eq:ConvexProp_eq}, if the expansion point $x_{\left(k-1\right)}$ is an admissible point to  problem~$P$ then, $x_{\left(k-1\right)}$ is also an admissible point to problem~$P_{\mathrm{cvx}}^{\left(k\right)}$, making the convexified problem~$P_{\mathrm{cvx}}^{\left(k\right)}$ feasible as $x_{\left(k-1\right)}\in\mathcal{D}_{P_{\mathrm{cvx}}^{\left(k\right)}}\Rightarrow\mathcal{D}_{P_{\mathrm{cvx}}^{\left(k\right)}}\neq\varnothing$\qed
\end{pf}

\begin{pf}[Theorem \ref{th:main_theorem}]
	The corollary of Proposition~\ref{prop:feasible} shows that if the first convexified problem~$P_{\mathrm{cvx}}^{\left(1\right)}$ is feasible then its solution, $x_{\left(1\right)}$, is an admissible point of the original non-convex problem~$P$. If that first solution, $x_{\left(1\right)}$, is used as the expansion point to obtain the next convexified problem~$P_{\mathrm{cvx}}^{\left(2\right)}$, Proposition~\ref{prop:admissible} ensures that this second problem is also feasible. If we continue to use the previous problem's solution as the expansion point to obtain the next convexified problem it follows that all the following convexified problems~$P_{\mathrm{cvx}}^{\left(k\right)}$ are feasible and all their solutions are admissible points of the original non-convex problem~$P$, $x_{\left(k\right)}\in\mathcal{D}_{P}$. Recursive feasibility is then established.
	
	Given Eq.~\eqref{eq:ConvexProp_geq} and \eqref{eq:ConvexProp_eq},
	\begin{equation}
	f_{0}\left(x_{\left(k\right)}\right) \leq f_{\text{cvx},0}^{\left(k-1\right)}\left(x_{\left(k\right)}\right)
	\leq f_{\text{cvx},0}^{\left(k-1\right)}\left(x_{\left(k-1\right)}\right)= f_{0}\left(x_{\left(k-1\right)}\right),
	\end{equation}
	and it is shown that the cost between iterations is non-strictly decreasing, \emph{i.e.}, making the SCP procedure an infinite descent algorithm.
	
	As the problem is lower bounded, $f_{0}\left(x\right)> -\infty$, the sequence must be convergent:
	\begin{equation}
	f_{\text{cvx},0}^{\left(\infty-1\right)}\left(x_{\left(\infty\right)}\right)=f_{\text{cvx},0}^{\left(\infty\right)}\left(x_{\left(\infty+1\right)}\right)=f_{0}\left(x_{\left(\infty\right)}\right),
	\end{equation}
	\begin{equation}
	x_{\left(\infty\right)}=x_{\left(\infty+1\right)}.
	\end{equation}
	As $x_{\left(\infty\right)}$ is the solution to problem~$P_{\mathrm{cvx}}^{\left(\infty\right)}$ it automatically satisfies the
	KKT conditions for problem~$P_{\mathrm{cvx}}^{\left(\infty\right)}$.
	Then, as $x_{\left(\infty\right)}$ is an admissible point of the non-convex problem~$P$ it satisfies its primal feasibility conditions, Eqs~\eqref{eq:primal_ineq} and \eqref{eq:primal_eq}. Using Eq.~\eqref{eq:ConvexProp_geq} and \eqref{eq:ConvexProp_eqDerivative} it is shown that
	\begin{equation}
	f_{\text{cvx}}^{\left(\infty\right)}\left(x_{\left(\infty\right)}\right)=f\left(x_{\left(\infty\right)}\right),
	\end{equation}
	\begin{equation}
	\nabla f_{\text{cvx}}^{\left(\infty\right)}\left(x_{\left(\infty\right)}\right)= \nabla f\left(x_{\left(\infty\right)}\right),
	\end{equation}
	and thus $x_{\left(\infty\right)}$, and its associated optimal dual variables $\nu_{\left(\infty\right)},\mu_{\left(\infty\right)}$, also meet the dual feasibility, (Eq.~\ref{eq:dual}) complementary slackness (Eq.~\ref{eq:compl}), and stationarity (Eq.~\ref{eq:stationarity}) conditions of the non-convex problem~$P$. This allows to conclude that $x_{\left(\infty\right)}$ is a locally optimal solution to the non-convex problem~$P$: $x_{\left(\infty\right)}=x^{\star}_{P}$\qed
\end{pf}
Two brief remarks are now offered about the proposed SCP. 1) to ensure that the first convexified problem in the sequence is feasible, an admissible initial guess, $x_{\left(0\right)}\in\mathcal{D}_{P}$, can be used to obtain the first convexified problem. If an admissible initial guess is not available, then the first problem in the SCP procedure may be infeasible, even if the original non-convex problem~$P$ is feasible. 2) the SCP procedure is local. The final solution potentially depends on the initial guess provided. Multiple initial guesses can be used to seed multiple SCP procedures, and eventually select the converged solution with the lowest cost.


\section{Approaches to obtain inner-convex approximations}
\label{sec:Methods}
In the literature we can find several methods to generate inner-convex approximations. These are briefly reviewed in this section.

If a function is expressed as the difference of two convex functions $c_{1}\left(x\right),c_{2}\left(x\right):\mathbb{R}^{n}\rightarrow\mathbb{R}$,
\begin{equation}
f\left(x\right)=c_{1}\left(x\right)-c_{2}\left(x\right).
\label{eq:dc}
\end{equation}
an inner-convex approximation can be obtained by linearizing the concave component, \emph{i.e.}, $-c_{2}\left(x\right)$, as follows:
\begin{equation}
f^{\left(x_{e}\right)}_{\text{cvx}}\left(x\right)=c_{1}\left(x\right)-c_{2}\left(x_{e}\right)-\delta x^{T}\nabla c_{2}\left(x_{e}\right).
\label{eq:cccp}
\end{equation}
These difference of convex functions are typically referred to as d.c., and the resulting SCP obtained when the concave component is linearized is known as the convex-concave procedure \cite{Yuille2003,Lipp2016}. This convex-concave procedure is thus a particular case of the proposed SCP with inner-convex approximations.

Problem formulations with d.c. functions naturally appear in many optimization problems \cite{Tuy1995,Tao2005}.  A common example in trajectory optimization are keep-out zones, which are usually formulated as concave inequality constraints (\emph{i.e.}, a d.c. with $c_{1}\left(x\right)=0$). Several authors have rediscovered the recursively feasible and convergent properties that an SCP enjoys when concave inequalities are linearized. \cite{Morgan2016,Liu2014,Virgili2018_IJRR}.

It has been shown that d.c. decompositions exist for any twice differentiable function $f\in\mathcal{C}^{2}$  \cite{Hartman1959,Hiriart1985}, and several approaches exists to find them.
For example, if $f$ is a polynomial, a d.c decomposition that is made by convex sum-of-squares polynomials (sos-convex in the following) can be found \cite{Ahmadi2017}. There are infinite sos-convex representations and thus in \cite{Ahmadi2017} different methods to optimize the sos-convex decomposition for a faster rate of convergence are proposed. Interestingly enough, finding these optimized sos-convex decompositions is a semidefinite programming problem (second order cone or linear programming problem with some additional relaxations) \cite{Ahmadi2017}.

Another approach to find d.c. decompositions is the regularization method \cite{Hiriart1985,Aleksandrov2012,Dinh2012,Wang2012b,Scutari2014}. In a regularization, a d.c. decomposition is obtained by using a sufficiently strong convex function $\phi\left(x\right):\mathbb{R}^{n}\rightarrow\mathbb{R}$ such that
\begin{equation}
f\left(x\right)=\underbrace{ \phi\left(x\right)}_{c_{1}\left(x\right)}-\underbrace{\left(\phi\left(x\right)-f\left(x\right)\right)}_{c_{2}\left(x\right)},\text{ for }x\in\mathcal{T},
\label{eq:regularization}
\end{equation}
is d.c in the trust region $\mathcal{T}$.

If the $f$ function is Lipschitz continuous on $x\in\mathcal{T}$, with $K$ denoting the Lipschitz constant, a d.c. representation of $f$ over $\mathcal{T}$ is 
\begin{equation}
f\left(x\right)=\underbrace{K\left\lVert x\right\rVert^{2}}_{c_{1}\left(x\right)}-\underbrace{\left(K\left\lVert x\right\rVert^{2}-f\left(x\right)\right)}_{c_{2}\left(x\right)},\text{ for }x\in\mathcal{T}.
\end{equation}

Using the same principle, and using the Lipschitz constant of the Hessian,  inner-convex approximations based on higher-order regularizations are also possible. For example, \cite{Nesterov2006} introduces a cubic regularization that enjoys a quadratic rate of convergence.

Up to this point, we have shown that inner-convex approximations for d.c. functions are easy to obtain. However, if the function is not d.c., either a semidefinite program needs to be solved to obtain a sos-convex decomposition or, for regularization-based approaches, we need to estimate the Lipschitz constant.



\section{Taylor-based inner-convex approximations}
\label{sec:IC_Taylor}
A new computationally inexpensive method to obtain an inner-convex approximation of a non-convex analytic function is here presented. This approach can be seen as a high order regularization based on the convexification of the high-order terms of a Taylor series expansion. With these inner-convex approximations the SCP's local rate of convergence is superlinear, or even quadratic when converging to a non-degenerate minima.

To simplify the notation of the high-order terms of the Taylor series let's use a multi-index notation, with an $n$-tuple of natural numbers $\mathbb{N}$, denoted by $\alpha$ \cite[Section~1.1]{Raymond1991}:
\begin{subequations}
\begin{equation}
\alpha=\left(\alpha_{1},\alpha_{2}\ldots \alpha_{n}\right),\quad \text{ with }\alpha_{j}\in\mathbb{N}
\end{equation}
\begin{equation}
\left\lvert\alpha\right\rvert=\alpha_{1}+\alpha_{2}+\cdots+\alpha_{n}
\end{equation}
\begin{equation}
\alpha!=\alpha_{1}!\alpha_{2}!\cdots\alpha_{n}!
\end{equation}
\begin{equation}
x^{\alpha}=x_{1}^{\alpha_{1}}x_{2}^{\alpha_{2}}\cdots x_{n}^{\alpha_{n}}
\end{equation}
\begin{equation}
\partial^{\alpha}f=\frac{\partial^{\left\lvert\alpha\right\rvert}f}{\partial x_{1}^{\alpha_{1}}\partial x_{2}^{\alpha_{2}}\cdots \partial x_{n}^{\alpha_{n}}}
\end{equation}
\end{subequations}
Using a multi-index notation, the Taylor series of the analytic function $f:\mathbb{R}^{n}\rightarrow\mathbb{R}$ is written as:
\begin{equation}
f\left(x\right)=\sum_{\left\lvert\alpha\right\rvert=0}^{\infty}{\frac{\partial^{\alpha}f\left(x_{e}\right)}{\alpha!}}\delta x^{\alpha}.
\end{equation}
To obtain an inner-convex approximation of $f$ we propose to convexify all of the terms corresponding to the derivatives of second order and higher. The Hessian, $\bm{H}$, as a symmetric matrix, can be decomposed in a positive and a negative semi-definite part \cite[Section~7.2]{Meyer2000}:
\begin{equation}
\bm{H}=\bm{H}^{+}+\bm{H}^{-}
\end{equation}
\begin{equation}
\bm{H}^{+}\succeq 0,\quad\bm{H}^{-}\preceq 0.
\end{equation}
The positive semi-definite part is recovered using the matrix of eigenvalues $\bm{\Lambda}$ and the matrix of eigenvectors $\bm{W}$ of the Hessian:
\begin{equation}
\bm{H}^{+}=\bm{W}\bm{\Lambda}^{+}\bm{W}^{T},
\end{equation}
where $\bm{\Lambda}^{+}$ is the matrix of eigenvalues containing only the positive ones. 

The higher order terms ($\left\lvert\alpha\right\rvert\ge 3$) of the Taylor series form homogeneous, polynomial forms
\begin{equation}
\sum_{\left\lvert\alpha\right\rvert}\frac{\partial^{\alpha}f\left(x_{e}\right)}{\alpha!}\delta x^{\alpha}=\bm{T}\delta x^{\alpha}=\sum_{j_{1}\ldots j_{\left\lvert\alpha\right\rvert}=1}^{n}{T_{j_{1}\ldots j_{\left\lvert\alpha\right\rvert}}\delta x_{j_{1}}\ldots \delta x_{j_{\left\lvert\alpha\right\rvert}}},
\end{equation}
with $\bm{T}\in\mathbb{R}^{n\times n\times \cdots \times n}$ denoting a supersymmetric tensor---invariant under permutation on any of its indices.

If the eigenvalues and eigenvectors of these tensors could be extracted a similar procedure to the one followed to convexify the Hessian could be pursued.
Although the eigenvalues of tensors can be computed using a sum-of-squares and a semidefinite  programming approach \cite{Cui2014} (showing a connection with sos-convex decompositions, Lasserre's relaxations, and non-convex programming), the computational cost quickly becomes prohibitive as the order, or number of variables and constraints grows.

Our proposed approach takes a computational efficient route and a convex overestimation of the terms of third order or higher, $\mathrm{cvx}\left(\bm{T}\delta x^{\alpha}\right)$, is obtained as follows:
\begin{equation}
\mathrm{cvx}\left(\bm{T}\delta x^{\alpha}\right)=\sum_{i=1}^{n}\left({\left\lvert T_{\mathrm{diag},i} \delta x_{i}^{\left\lvert\alpha\right\rvert}\right\rvert^{+}+{T_{\mathrm{cvx},i}\left\lvert \delta x_{i}^{\left\lvert\alpha\right\rvert}\right\rvert}}\right).
\end{equation}
The $T_{\mathrm{diag},i}$ coefficient is the element of the tensor's main diagonal associated with $\delta x_{i}$,
\begin{equation}
T_{\mathrm{diag},i}=T_{i,i\ldots i}.
\end{equation}
As $\left\lvert \cdot \right\rvert^{+}$ denotes the positive part operator, the expression ${\left\lvert T_{\mathrm{diag},i} \delta x_{i}^{\left\lvert\alpha\right\rvert}\right\rvert^{+}}$ is convex and it is an overestimator of the $\bm{T}\delta x^{\alpha}$ terms associated with the $\bm{T}$ tensor's main diagonal.

The $T_{\mathrm{cvx},i}$ coefficient is computed as the sum of the absolute values of all of the coefficients in $\bm{T}$ associated with $\delta x_{i}$ (with the exception of the coefficients in the main diagonal):
\begin{equation}
T_{\mathrm{cvx},i}{=}\sum_{j_{1}\ldots j_{\left\lvert\alpha\right\rvert}=1}^{n}
\begin{cases}
0&\text{if all }j_{1}\ldots j_{\left\lvert\alpha\right\rvert}\text{ are different than }i\\
0&\text{if all }j_{1}\ldots j_{\left\lvert\alpha\right\rvert}\text{ are equal to }i\\
{\left\lvert T_{j_{1}\ldots j_{\left\lvert\alpha\right\rvert}}\right\rvert}&\text{otherwise}
\end{cases}
\end{equation}
As $T_{\mathrm{cvx},i}\geq 0$ the term $T_{\mathrm{cvx},i}\left\lvert \delta x_{i}^{\left\lvert\alpha\right\rvert}\right\rvert$ is convex and it follows that
\begin{equation}
\mathrm{cvx}\left(\bm{T}\delta x^{\alpha}\right)\geq\bm{T}\delta x^{\alpha}.
\end{equation}
If the coefficients $T_{\mathrm{cvx},i}$ and $T_{\mathrm{diag},i}$ of a particular order, $\left\lvert\alpha\right\rvert$, are denoted by $T^{\left\{\left\lvert\alpha\right\rvert \right\}}_{\mathrm{cvx},i}$ and $T^{\left\{\left\lvert\alpha\right\rvert \right\}}_{\mathrm{diag},i}$, the proposed inner-convex approximation of $f$, satisfying the conditions in Eq.~\eqref{eq:ConvexProp}, is
\begin{equation}
\begin{split}
f^{\left(x_{e}\right)}_{\text{cvx}}\left(x\right)&=f\left(x_{e}\right)+\delta x^{T}\nabla f\left(x_{e}\right)+\frac{1}{2}\delta x^{T}\bm{H}^{+}\left(x_{e}\right)\delta x\\
&+\sum_{\left\lvert\alpha\right\rvert=3}^{\left\lvert\alpha\right\rvert=\infty}{\sum_{i=1}^{n}\left({\left\lvert T^{\left\{\left\lvert\alpha\right\rvert \right\}}_{\mathrm{diag},i} \delta x_{i}^{\left\lvert\alpha\right\rvert}\right\rvert^{+}}+{T^{\left\{\left\lvert\alpha\right\rvert \right\}}_{\mathrm{cvx},i}\left\lvert \delta x_{i}^{\left\lvert\alpha\right\rvert}\right\rvert}\right)}.
\label{eq:gcvx}
\end{split}
\end{equation}

\subsection{Truncating the series}\label{sec:trust_regions}
The inclusion of high order terms might need to be truncated as the Taylor series may be infinite or the computational effort required to compute high-order terms may be too high. 
To truncate the series and preserve the inner-convex properties a high-order regularization can be added as follows,
\begin{equation}
f^{\left(x_{e}\right)}_{\text{cvx}}\left(x,M\right)=\overbrace{f^{\left(x_{e}\right)}_{\text{cvx}}\left(x\right)}^{\text{truncated at order } \left\lvert\alpha\right\rvert}+\frac{M}{\left(\left\lvert\alpha\right\rvert+1\right)!}\left\lVert\delta x\right\rVert^{\left\lvert\alpha\right\rvert+1},
\end{equation}
for $M\geq0$. 

If the condition $f_{\mathrm{cvx}}^{\left(k-1\right)}\left(x_{\left(k\right)},M\right)\geq f\left(x_{\left(k\right)}\right)$ is not met after solving the convex problem $P_{\mathrm{cvx}}^{\left(k\right)}\left(M\right)$, the problem is re-solved with a new, larger $M$, \emph{i.e.}, $M_{\mathrm{new}}>M_{\mathrm{old}}$. As long as the  $\left\lvert\alpha\right\rvert+1$ order derivative of $f$ is Lipschitz continuous the $f_{\mathrm{cvx}}^{\left(k-1\right)}\left(x_{\left(k\right)},M\right)\geq f\left(x_{\left(k\right)}\right)$ condition will be eventually met by a large-enough $M$. Adding the high-order regularization may potentially reduce the rate of convergence and therefore it is advised to start with $M=0$, only increasing $M$ when strictly required. 

This truncation with a high-order regularization can be seen as an extension of the cubic regularization approach developed by Nesterov \cite{Nesterov2006}.

\subsection{Quadratic rate of convergence}
The proposed inner-convex approximation of the Taylor series provides the SCP with a local superlinear rate of converge, which can be proven to be quadratic when converging to a non-degenerate minima.

It has been shown that a linearization-based SCP is superlinearly convergent around its convergence point $x_{\left(\infty\right)}$ \cite{Mao2018}. As the SCP converges,
\begin{equation}
\lim\limits_{k\rightarrow\infty}\delta x=0,
\end{equation}
making the high-order terms become increasingly small
\begin{equation}
f^{\left(x_{e}\right)}_{\text{cvx}}\left(x\right)=f\left(x_{e}\right)+\delta x^{T}\nabla f\left(x_{e}\right)
+O\left(\left\lVert\delta x\right\rVert^{2}\right),
\end{equation}
and thus the Taylor-based convexification inherits the superlinear convergence of linearization-based SCPs.

But as the Taylor-based convexification preserves the Hessian, as long as $\bm{H}\left(x\right)\succ0$, the resulting SCP exhibits a quadratic rate of convergence for non-degenerate minima.
\begin{thm}
	The SCP exhibits a quadratic rate of convergence,
	\begin{equation*}
	\left\lVert x_{\left(k+1\right)}-x_{\left(k\right)}\right\rVert\leq
	\gamma_{\left(k\right)}\left\lVert x_{\left(k\right)}-x_{\left(k-1\right)}\right\rVert^{2},
	\end{equation*}
	for some $\gamma_{\left(k\right)}>0$, around a non-degenerate local minima, $\bm{H}\left(x\right)\succ0$, of Problem $P$.
\end{thm}
The cubic regularization approach proposed by Nesterov \cite{Nesterov2006} also has a quadratic rate of convergence, and the proof provided here is based on the proof provided for that case. 
\begin{pf} 
Let us start by introducing some additional notation.
Let the spectrum of a of an $n\times n$ symmetric matrix $\bm{A}$ be denoted by $\left\{\lambda_{i}\left(\bm{H}\right)\right\}_{i=1}^{n}$. Also, let's assume that the eigenvalues are ordered in a decreasing sequence:
\begin{equation}
\lambda_{1}\left(\bm{A}\right)\geq\ldots\geq\lambda_{n}\left(\bm{A}\right),
\end{equation}
and that the norm of a real symmetric matrix is defined with its standard spectral form as \cite[Chapter~5]{Horn1990}
\begin{equation}
\left\lVert\bm{A}\right\rVert=\sqrt{\lambda_{1}\left(\bm{A}^{T}\bm{A}\right)}.
\end{equation}
Having introduced this notation, let us suppose that the optimal dual variables at each iteration of the SCP are known, then, the Lagrangian of the problem with the non-convex functions is
\begin{equation}
L_{\left(k\right)}\left(x\right)=f_{0}\left(x\right)+\sum_{i=1}^{p}{\nu_{\left(k\right)i}^{\star} f_{i}\left(x\right)}+\sum_{i=1}^{q}{\mu_{\left(k\right),i}^{\star} h_{i}\left(x\right)}.
\label{eq:Lagrangian}
\end{equation}
Solving the constrained problem $P$ is equivalent to solving an unconstrained problem with the Lagrangian in Eq.~\ref{eq:Lagrangian} as the cost function (dual problem). Without any loss of generality, let's use this unconstrained version of the problem for the proof.

Further assume that the Hessian of the Lagrangian, $\bm{H}_{L}$ is locally Lipschitz continuous,
\begin{equation}
\left\lVert \bm{H}_{L}\left(x_{2}\right)-\bm{H}_{L}\left(x_{1}\right)\right\rVert\leq K_{L} \left\lVert x_{2}-x_{1}\right\rVert.
\end{equation}
Integrating both sides of the inequality yields:
\begin{equation}
\left\lVert \nabla L\left(x_{2}\right)-\nabla L\left(x_{1}\right)-\bm{H}_{L}\left(x_{1}\right)\left(x_{2}-x_{1}\right)\right\rVert
\leq \frac{K_{L}}{2} \left\lVert x_{2}-x_{1}\right\rVert^{2}.
\label{eq:L_int}
\end{equation}
Using the inner-convex approximations at iteration $k+1$ the convexified Lagrangian is,
\begin{equation}
L^{\left(k\right)}_{\mathrm{cvx}}\left(x\right)=f^{\left(k\right)}_{\mathrm{cvx},0}\left(x\right)+\sum_{i=1}^{p}{\nu_{\left(k+1\right),i}^{\star} f^{\left(k\right)}_{\mathrm{cvx},i}\left(x\right)}
+\sum_{i=1}^{q}{\mu_{\left(k+1\right),i}^{\star} h_{i}\left(x\right)},
\end{equation}
and, as a sum of convex and affine functions, $L^{\left(k\right)}_{\mathrm{cvx}}$ is convex and takes the form of a convexified Taylor series as defined in Eq.~\eqref{eq:gcvx}.

The stationary condition (see Eq.~\eqref{eq:stationarity}) imposes that the gradient of the convexified Lagrangian $L_{\mathrm{cvx}}^{\left(k\right)}$ at the optimal point $x_{\left(k+1\right)}$ is zero. 
\begin{equation}
\nabla L_{\mathrm{cvx}}^{\left(k\right)}\left(x_{\left(k+1\right)}\right)=0.
\end{equation}
Expanding it as Taylor series yields
\begin{equation}
\nabla L\left(x_{\left(k\right)}\right)+\bm{H}^{+}_{L}\left(x_{\left(k\right)}\right)\left(x_{\left(k+1\right)}-x_{\left(k\right)}\right)+\ldots=0.
\end{equation}
If it is further assumed that the current iterate $x_{\left(k\right)}$ is close to a non-degenerate local minima with $\bm{H}_{L}\succ 0$. With this assumption, and noting that the high-order terms of $L_{\mathrm{cvx}}^{\left(k\right)}$ are all convex, it follows that 
\begin{equation}
\left\lVert x_{\left(k+1\right)}-x_{\left(k\right)}\right\rVert=\left\lVert \left(\bm{H}_{L}\left(x_{\left(k\right)}\right)+\ldots\right)^{-1}\nabla L\left(x_{\left(k\right)}\right)\right\rVert
\leq
\frac{\left\lVert \nabla L\left(x_{\left(k\right)}\right)\right\rVert}{\lambda_{n}\left(\bm{H}_{L}\left(x_{\left(k\right)}\right)\right)}.
\label{eq:optimality}
\end{equation}
Additionally, the Hessian of the convexified Lagrangian $\bm{H}_{L^{\left(k\right)}_{\mathrm{cvx}}}$ is locally Lipschitz continuous
\begin{equation}
\left\lVert \bm{H}_{L^{\left(k\right)}_{\mathrm{cvx}}}\left(x\right)-\bm{H}_{L^{\left(k\right)}_{\mathrm{cvx}}}\left(x_{\left(k\right)}\right)\right\rVert\leq K_{L^{\left(k\right)}_{\mathrm{cvx}}} \left\lVert x-x_{\left(k\right)}\right\rVert,
\end{equation}
allowing to bound the high order terms
\begin{multline}
L_{\mathrm{cvx}}^{\left(k\right)}\left(x\right)\leq
L\left(x_{\left(k\right)}\right)+\nabla L\left(x_{\left(k\right)}\right)\left(x-x_{\left(k\right)}\right)\\+
\left(x-x_{\left(k\right)}\right)^{T}\bm{H}_{L}\left(x_{\left(k\right)}\right)\left(x-x_{\left(k\right)}\right)+\frac{K_{L^{\left(k\right)}_{\mathrm{cvx}}}}{6}\left\lVert x-x_{\left(k\right)}\right\rVert^{3}.
\end{multline}
With this observation, the optimality condition in Eq.~\eqref{eq:optimality} yields the following inequality:
\begin{equation}
\left\lVert \nabla L\left(x_{\left(k\right)}\right)+ \bm{H}_{L}\left(x_{\left(k\right)}\right)\left(x_{\left(k+1\right)}-x_{\left(k\right)}\right)\right\rVert\leq \frac{1}{2}K_{L^{\left(k\right)}_{\mathrm{cvx}}}\left\lVert x_{\left(k+1\right)}-x_{\left(k\right)}\right\rVert^{2}.
\label{eq:traingle_2}
\end{equation}
Combining Eq.~\eqref{eq:L_int} with Eq.~\eqref{eq:traingle_2} and observing the triangle inequality yields the following expression:
\begin{equation}
\left\lVert \nabla L\left(x_{\left(k+1\right)}\right)\right\rVert\leq \frac{K_{L}+K_{L^{\left(k\right)}_{\mathrm{cvx}}}}{2}\left\lVert x_{\left(k+1\right)}-x_{\left(k\right)}\right\rVert^{2}.
\end{equation}
Therefore it follows that
\begin{equation}
\left\lVert x_{\left(k+1\right)}-x_{\left(k\right)}\right\rVert\leq\frac{\left\lVert \nabla L\left(x_{\left(k\right)}\right)\right\rVert}{\lambda_{n}\left(\bm{H}_{L}\left(x_{\left(k\right)}\right)\right)}
\leq
\frac{\left(K_{L}+K_{L^{\left(k\right)}_{\mathrm{cvx}}}\right)}{2\lambda_{n}\left(\bm{H}_{L}\left(x_{\left(k\right)}\right)\right)}\left\lVert x_{\left(k\right)}-x_{\left(k-1\right)}\right\rVert^{2},
\end{equation}
thus establishing the algorithm's local quadratic rate of convergence\qed
\end{pf}

\section{Application to a trajectory optimization problem}
\label{sec:Simulation}
To illustrate the use of the proposed SCP procedure with Taylor-based inner-convex approximations, an illustrative example is offered. In this example, inspired by the one provided in \cite{Mao2017b,Mao2018}, the trajectory of an aerial  vehicle subject to drag, control bounds, and a keep-out zone constraint is optimized. The dynamics of the system is written as
\begin{equation}
\ddot{r}=\frac{1}{m}\left(F-k_{d}\left\lVert\dot{r}\right\rVert\dot{r}\right)
\label{eq:dynamics},
\end{equation}
with $r\in\mathbb{R}^{3}$ denoting the position of the vehicle, $m$ its mass, $F\in\mathbb{R}^{3}$ its thrust, and $k_{d}$ the vehicle's drag coefficient.

The vehicle is subject to the following control constraint,
\begin{equation}
\left\lVert F\right\rVert\leq F_{\mathrm{max}},
\label{eq:Tmax}
\end{equation}
and it must respect the following keep-out zone constraint,
\begin{equation}
(r_{1}^2+r_{2}^2)^2+r_{3}^4-b^4-10r_{3}(r_{1}^2r_{2}-r_{2}^2r_{1})\geq 0.
\label{eq:KO}
\end{equation}
The initial and final conditions, as well as the time to complete the maneuver, $t_{f}$, are specified:
\begin{subequations}
\begin{equation}
r\left(0\right)=r_{0},\quad\dot{r}\left(0\right)=\dot{r}_{0},
\label{eq:init_state}
\end{equation}
\begin{equation}
r\left(t_{f}\right)=r_{f},\quad\dot{r}\left(t_{f}\right)=\dot{r}_{f}.
\label{eq:final_state}
\end{equation}
\end{subequations}
The cost to be minimized is expressed as
\begin{equation}
J=\int_{0}^{t_{f}}\left\lVert F\right\rVert dt.
\label{eq:cost}
\end{equation}
This optimal control problem is discretized with a direct transcription method with $N$ nodes, applying a first-order-hold on the accelerations $\ddot{r}$, and estimating the cost with a trapezoidal integration scheme. The resulting non-convex programming problem is:
\begin{align}
\text{minimize: } & J=\sum_{i=1}^{N-1}\frac{\left\lVert F^{[i+1]}\right\rVert-\left\lVert F^{[i]}\right\rVert}{2} \left(t^{[i+1]}-t^{[i]}\right)\label{eq:cost_simulation}\\
\text{subject to: } &  F^{[i]}=m\ddot{r}^{[i]}+k_{d}\left\lVert\dot{r}^{[i]}\right\rVert\dot{r}^{[i]}, \tag{\ref{eq:dynamics}}\\
& \left\lVert F^{[i]}\right\rVert\leq F_{\mathrm{max}}, \tag{\ref{eq:Tmax}}\\
\begin{split}
(r_{1}^{[i]2}{+}r_{2}^{[i]2})^2&{-}10r_{3}^{[i]}(r_{1}^{[i]2}r_{2}^{[i]}{-}r_{2}^{[i]2}r_{1}^{[i]}){+}r_{3}^{[i]4}{-}b^4{\geq} 0
\end{split} \tag{\ref{eq:KO}}\\
&r^{[1]}=r_{0},\quad\dot{r}^{[1]}=\dot{r}_{0} \tag{\ref{eq:init_state}}\\
&r^{[N]}=r_{f},\quad\dot{r}^{[N]}=\dot{r}_{f} \tag{\ref{eq:final_state}}
\end{align}
Note that this problem is of the class $P$ with $x=\left\lbrace\ddot{r}^{\left[1\right]}\ldots\ddot{r}^{\left[N-1\right]}\right\rbrace$. The cost as well as the control and keep-out constraints are non-convex functions. 

For the cost and control constraints, the underlying non-convex functions have a Taylor series expansion with infinite terms and an inner-convex approximation truncated to order $\left\lvert\alpha\right\rvert =3$ is used during the SCP. No regularization is added.

The keep-out zone constraint can be re-written as follows
\begin{equation}
\underbrace{-(r_{1}^2+r_{2}^2)^2-r_{3}^4+b^4}_{\text{sos-concave}}\underbrace{+10r_{3}(r_{1}^2r_{2}-r_{2}^2r_{1})}_{\text{non-convex}}\leq 0.\label{eq:koz_d.c.}
\end{equation}
showing that it has a concave part---specifically sos-concave---and a non-convex part.
The Taylor series of the non-convex part is finite, reaching fourth order. To obtain inner-convex approximation of the non-convex part we convexify all the terms of the Taylor series expansion up to fourth order. Once the non-convex part of Eq.~\eqref{eq:koz_d.c.} has been convexified the function becomes d.c. and thus we know that the concave part needs to be linearized to obtain an inner-convex approximation.
This example problem shows how different methods can be used to generate inner-convex approximations.

A Monte Carlo analysis with a total of 1000 cases has been conducted. The parameters used are provided in Table~\ref{tab:parameters}.
\begin{table}
	\caption{\label{tab:parameters} Simulation parameters (non-dimensional).}
	\centering
	\begin{tabular}{cc}
		Parameter & Value \\
		\hline \hline
		\multicolumn{2}{c}{\emph{Problem}} \\
		Mass $m$ & 1 \\
		Final time $tf$ & 15 \\
		Drag coefficient $k_{d}$ & 0.25 \\
		Max. thrust $F_{\mathrm{max}}$ & 1.5	 \\
		Keep-out zone parameter $b$ & 	3.5\\
		Number of nodes $N$ & 25\\
		SCP stopping criterion $\epsilon$ & $0.01 f_{0}\left(\bm{x}^{\star}_{\left(k\right)}\right)$ \\
		Max. num. of SCP iterations  & 50\\
		\multicolumn{2}{c}{\emph{Monte Carlo}} \\
		Number of cases & 1000 \\
		Initial position $r_{0}$ & Random with $\left\lVert r_{0}\right\rVert=6$ \\
		Final position $r_{f}$ & $r_{f}=-r_{0}$ \\
		Initial velocity $\dot{r}_{0}$ & Random with $\left\lVert\dot{r}_{0}\right\rVert=1$ \\
		Final velocity $\dot{r}_{f}$ & Random with $\left\lVert\dot{r}_{f}\right\rVert=1$ \\
		\hline \hline
	\end{tabular}
\end{table}
Based on the randomized initial and final conditions an initial guess is analytically generated for each case. In particular, a bang-bang type control keeping the acceleration at two different constant values during the first and second halves of the maneuver is used to generate the initial guess. This solution meets the initial and final conditions but it may violate the thrust and/or keep-out zone constraints. As a consequence, the initial guess may not be an admissible solution of the non-convex problem. Therefore, we cannot guarantee a recursively feasible SCP procedure based on an inadmissible initial guess.

To generate an admissible solution from the initial guess a penalty method is used. This step serves to illustrate how inner-convex approximations can also be applied in conjunction with penalty methods and help obtain an admissible solution. In the penalty method version of the problem, the maximum thrust constraint is relaxed with a slack variable $s$,
\begin{subequations}
\begin{equation}
\left\lVert F^{[i]}\right\rVert- F_{\mathrm{max}}\leq s^{[i]},
\end{equation}
\begin{equation}
s^{[i]}\geq 0.
\end{equation}
\end{subequations}
and the cost function in Eq.~\eqref{eq:cost_simulation} replaced by
\begin{equation}
J=\sum_{n=1}^{N}s^{[i]}.
\end{equation}
The penalty method version of the problem is then solved using the same SCP approach until an admissible solution is obtained. Given that inner-convex approximations are used and that the problem's feasibility is guaranteed by the introduction of the slack variable, the sequence is bound to converge. However, the converged solution of the penalty version of the problem may be an inadmissible point of the original problem. 

When an admissible solution is found, the slack variable is removed, the original cost function in Eq.~\eqref{eq:cost_simulation} restored, and the SCP continues. As an admissible solution is now used as a seed, the remaining iterations of the SCP procedure are recursively feasible and guaranteed to converge to a locally optimal solution. 

\subsection{Illustrative result of one case}
\begin{figure*}
	\subfloat[]{\includegraphics[width=0.5\textwidth]{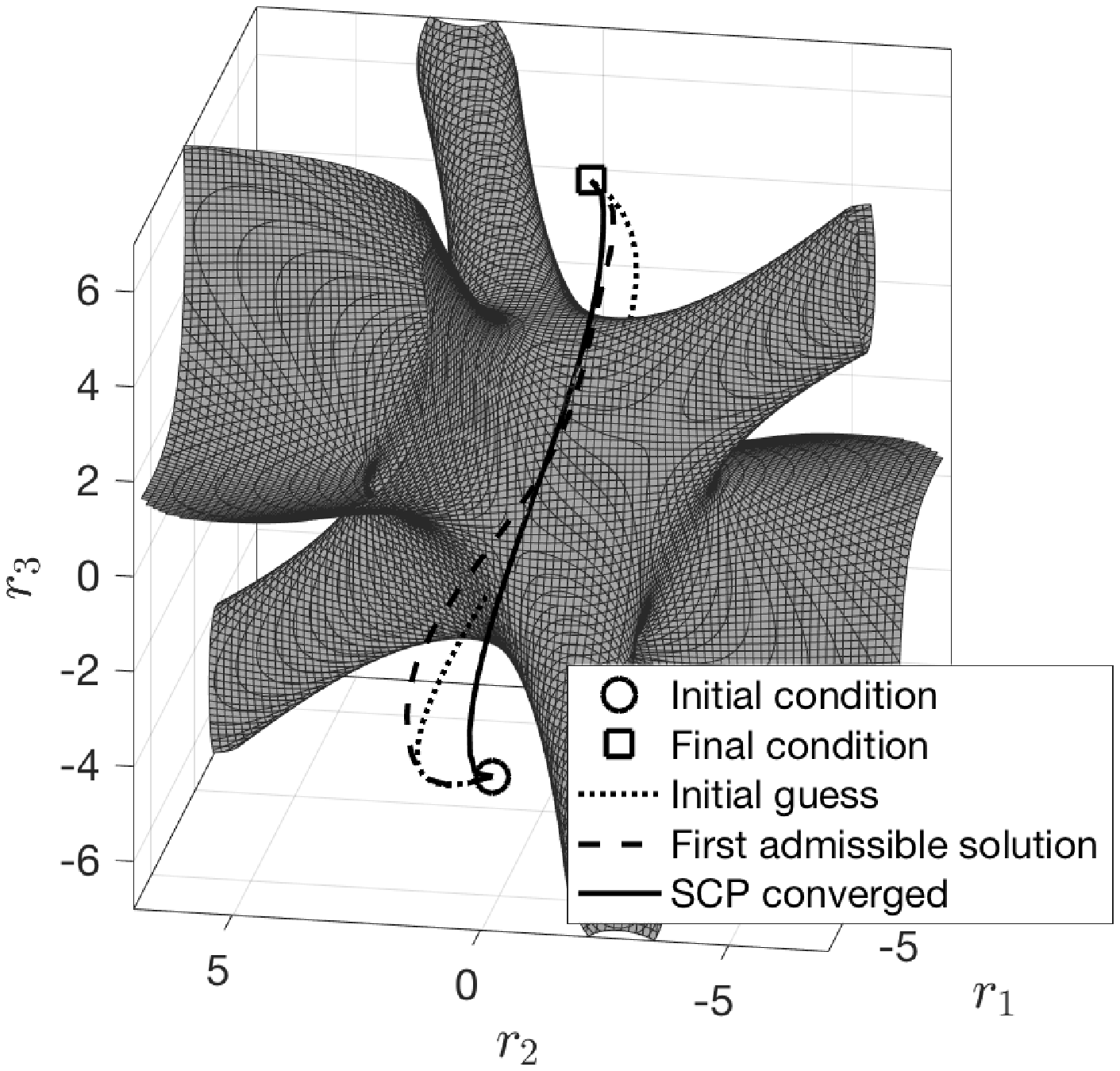}}
	\subfloat[]{\includegraphics[width=0.5\textwidth]{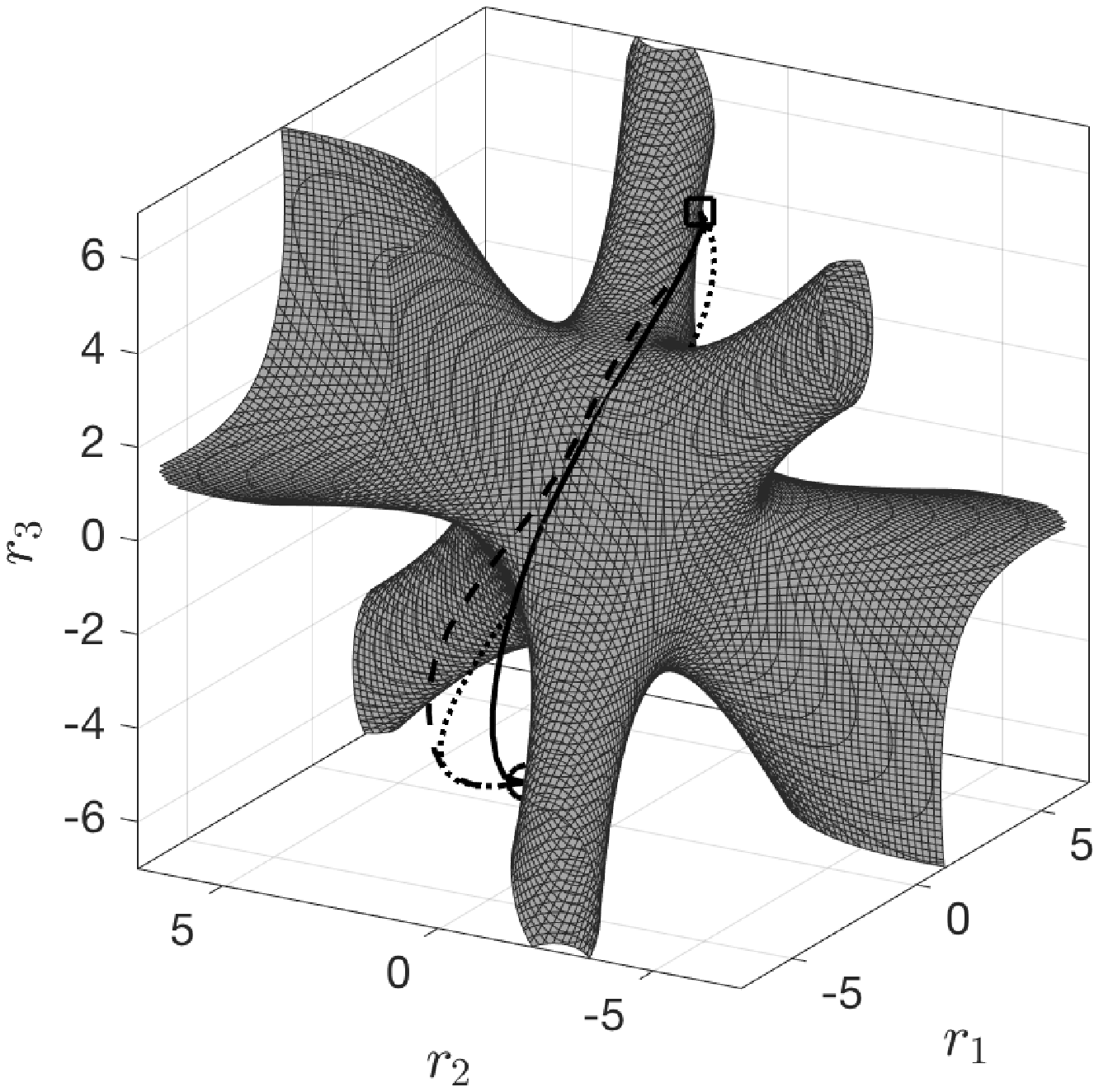}}\\
	\subfloat[]{\includegraphics[width=0.5\textwidth]{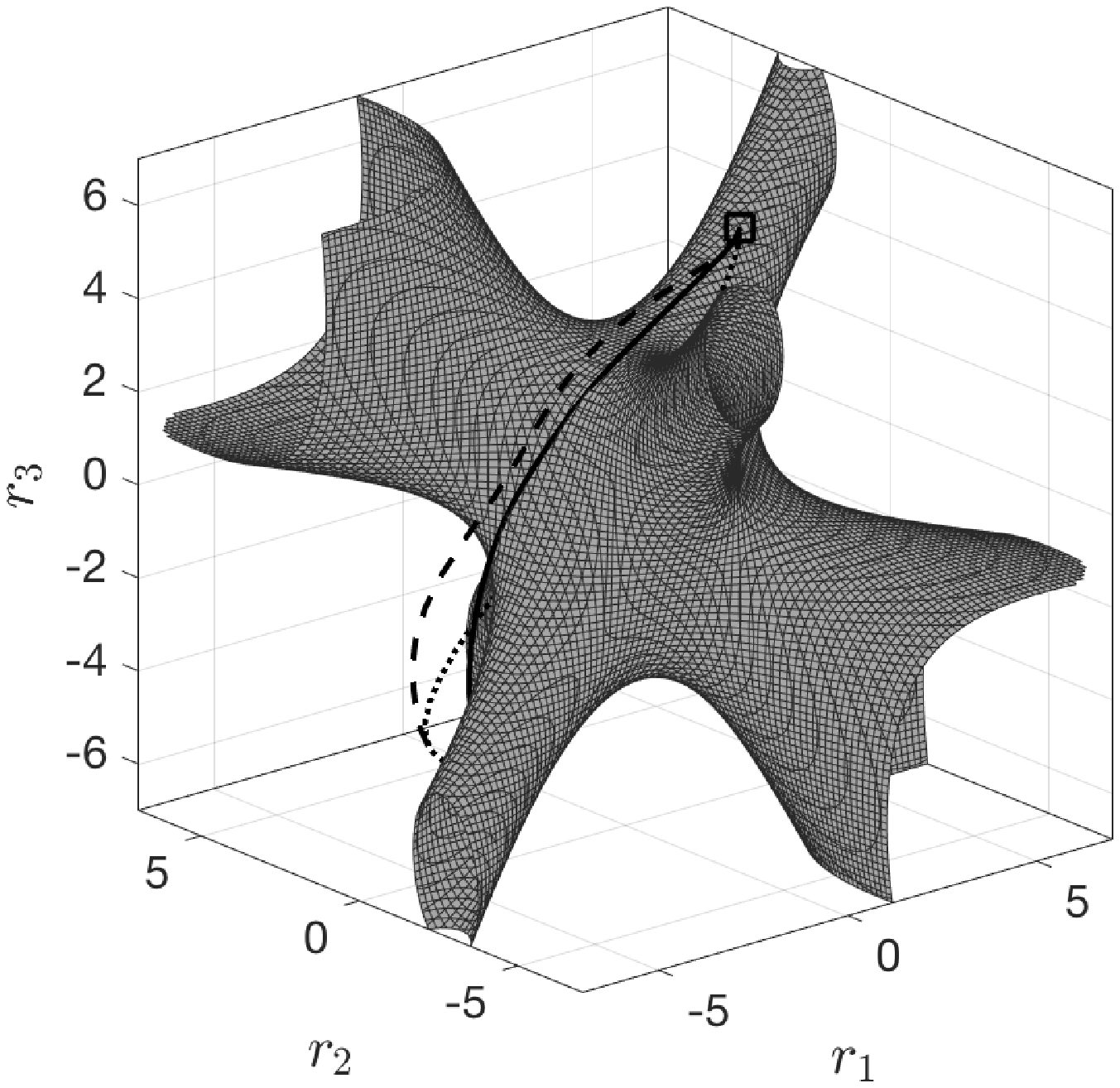}}
	\subfloat[]{\includegraphics[width=0.5\textwidth]{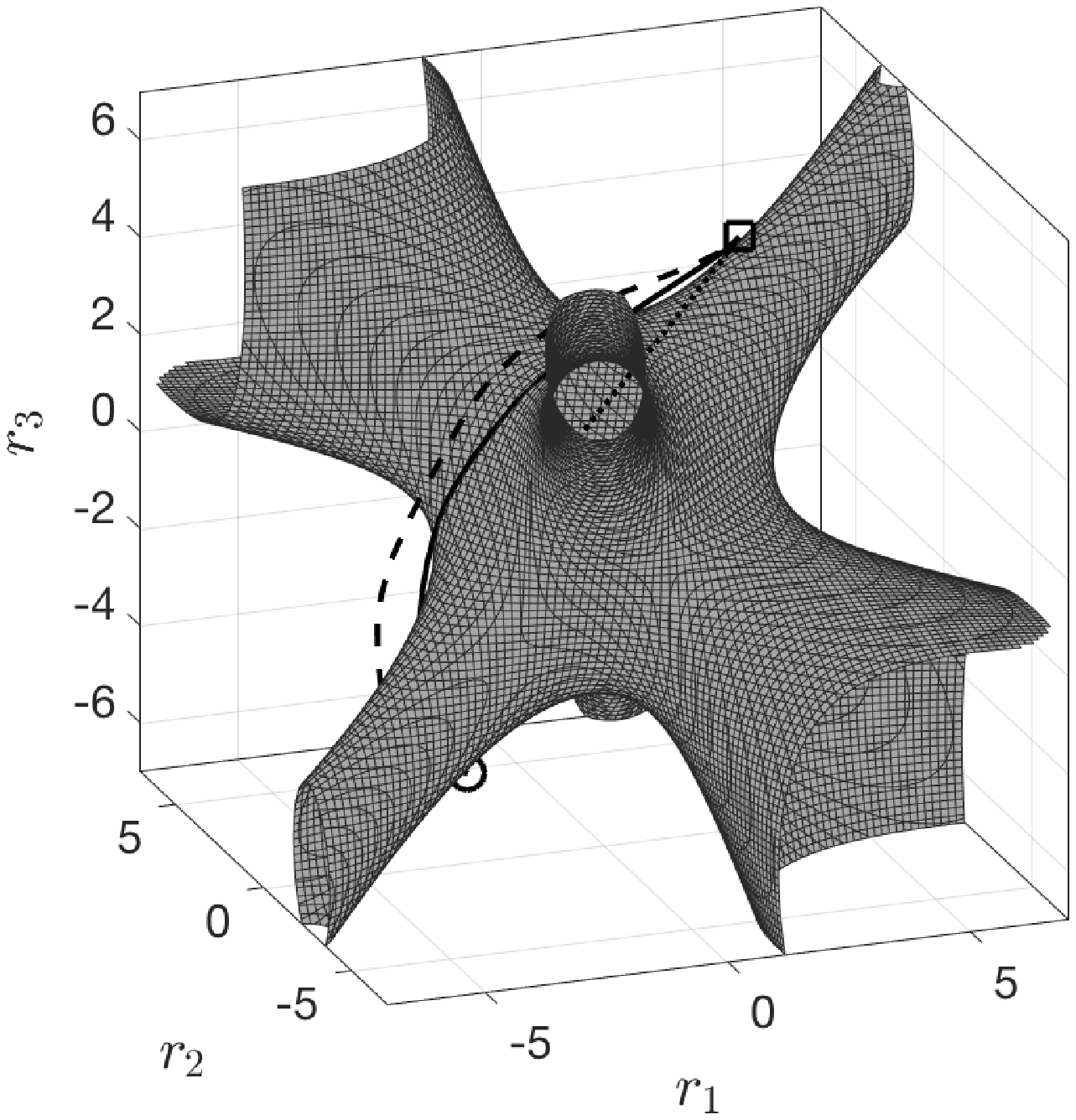}}
	\caption{Initial guess, first admissible solution, and converged trajectory at different viewing angles.\label{fig:trajectory}}
\end{figure*}

\begin{figure}
	\includegraphics[width=0.5\textwidth]{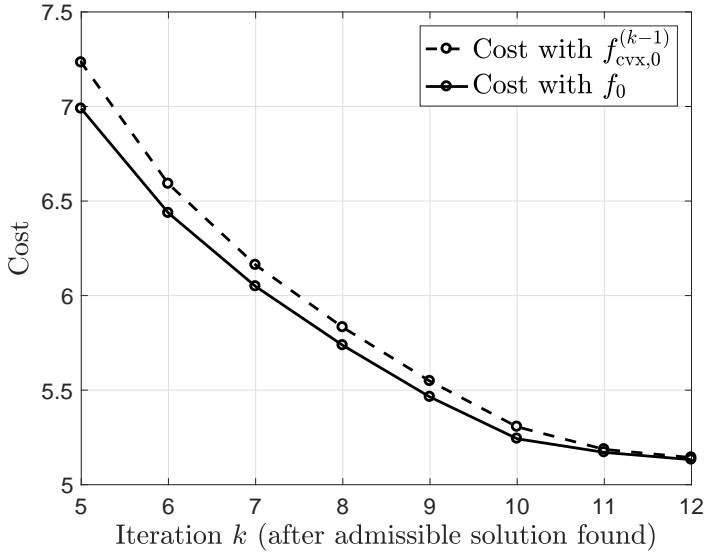}
	\caption{Evolution of the cost during the SCP.\label{fig:cost_evolution}}
\end{figure}

\begin{figure*}
	\subfloat[Thrust norm $\left\lVert F\right\rVert$.\label{fig:thrust}]{\includegraphics[width=0.5\textwidth]{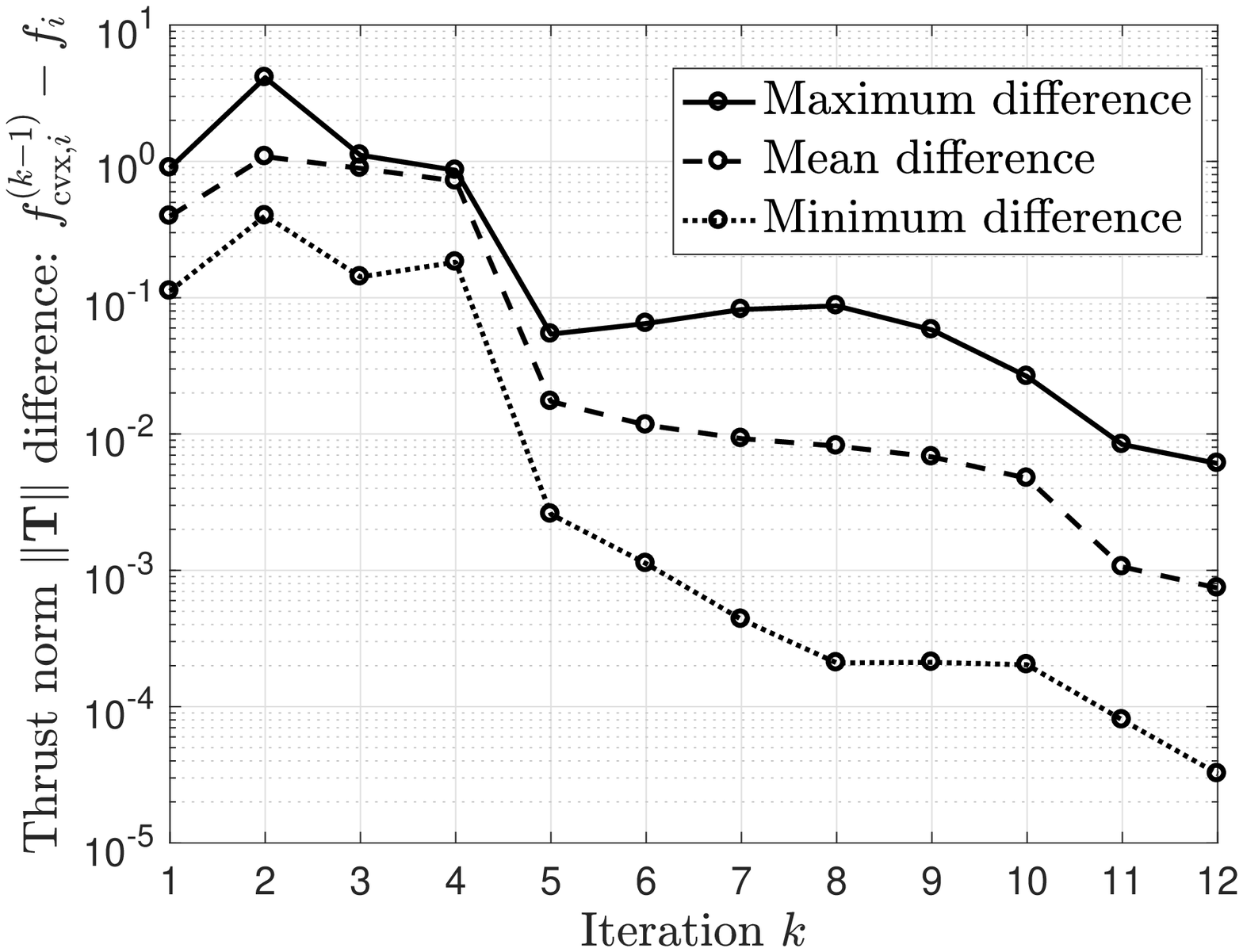}}	
	\subfloat[Keep-out zone constraint in Eq.~\eqref{eq:KO}.\label{fig:KO}]{\includegraphics[width=0.5\textwidth]{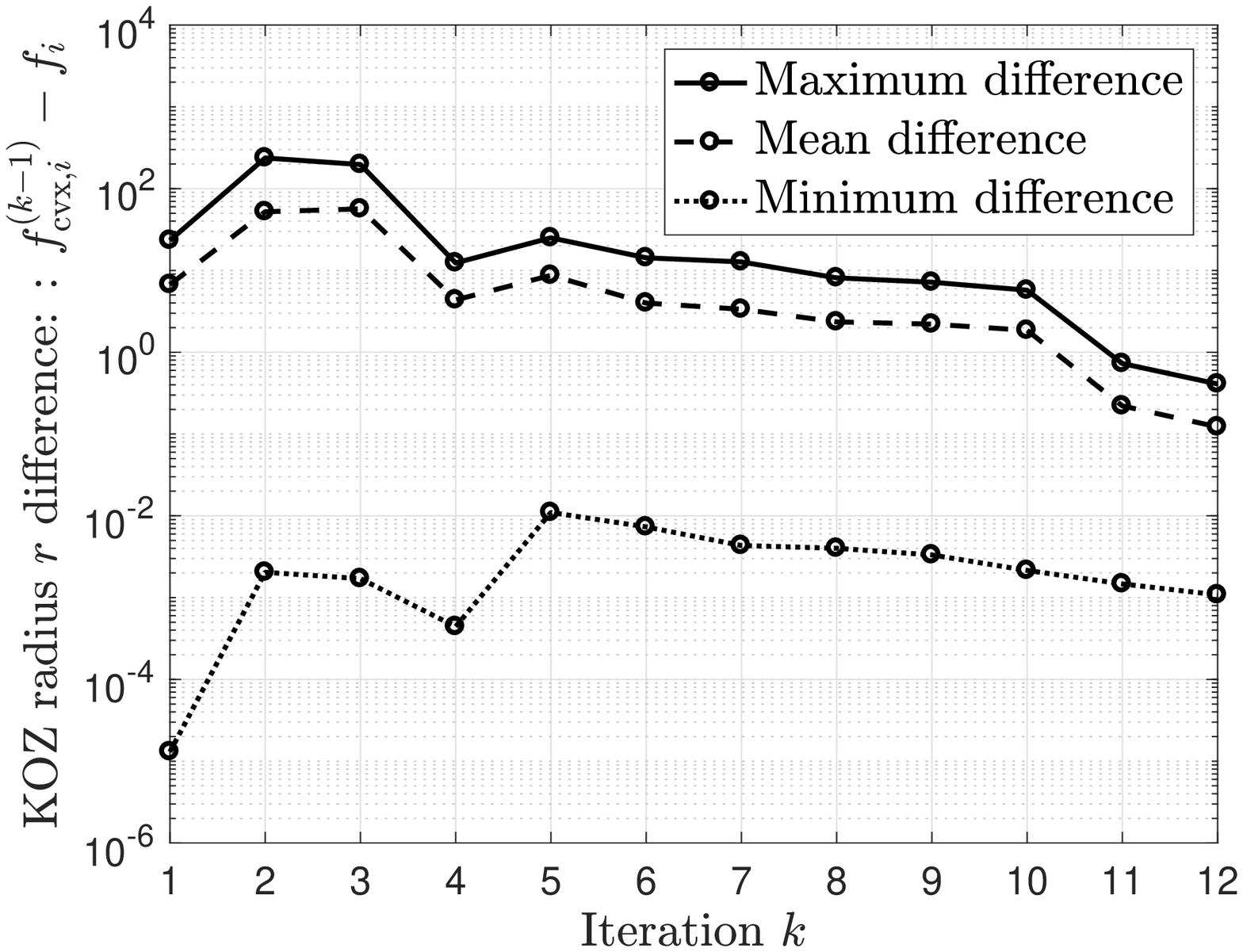}}\\
	\caption{Evolution of the constraint differences, $f_{\mathrm{cvx},i}^{(k-1)}-f_{i}\geq 0$, during the SCP.\label{fig:constraint_evolution}}
\end{figure*}

For illustration purposes, the results of one case with the initial and terminal conditions indicated in Table \ref{tab:parameters_one} are shown in Figs.~\ref{fig:trajectory}-\ref{fig:constraint_evolution}.
\begin{table}
	\caption{\label{tab:parameters_one} Initial and terminal condition of illustrative case.}
	\centering
	\begin{tabular}{cc}
		Parameter & Value \\
		\hline \hline
		Initial position $r_{0}$ &  $r_{0}=\left[-2.61,\, 0.53,\, -5.38\right]^{T}$ \\
		Final position $r_{f}$ & $r_{f}=-r_{0}$ \\
		Initial velocity $\dot{r}_{0}$ & $\dot{r}_{0}=\left[-0.62,\, 0.77,\, -0.14\right]^{T}$ \\
		Final velocity $\dot{r}_{f}$ & $\dot{r}_{f}=\left[0.64,\,
		0.75,\,
		0.15\right]^{T}$ \\
		\hline \hline
	\end{tabular}
\end{table}
The keep-out zone and the trajectories corresponding to the initial guess, the first admissible solution, and the converged solution are shown in Fig.~\ref{fig:trajectory}. In this particular case, four iterations with the penalty method were needed to produce an admissible solution. Eight more iterations were then required to converge down to the specified level, bringing the total number of convex problems solved to 12. Figure~\ref{fig:cost_evolution} shows the decreasing cost along with the inner-convex approximation of the cost, which is always overestimating it. The non-convex thrust and keep-out zone constraints are also overestimated by the inner-convex approximation, $f_{\mathrm{cvx},i}^{(k-1)}-f_{i}\geq 0$, as shown in Figs.~\ref{fig:thrust} and \ref{fig:KO}. In all of the SCP iterations within the 1000 cases, the thrust was always overestimated, showing that when the Taylor series is truncated a high order regularization may not be required in all cases.

\subsection{Results of the Monte Carlo analysis}
Analyzing the Monte Carlo results provides further insight into the proposed method, specifically with respect to its computational properties.

\begin{figure}
	\includegraphics[width=0.5\textwidth]{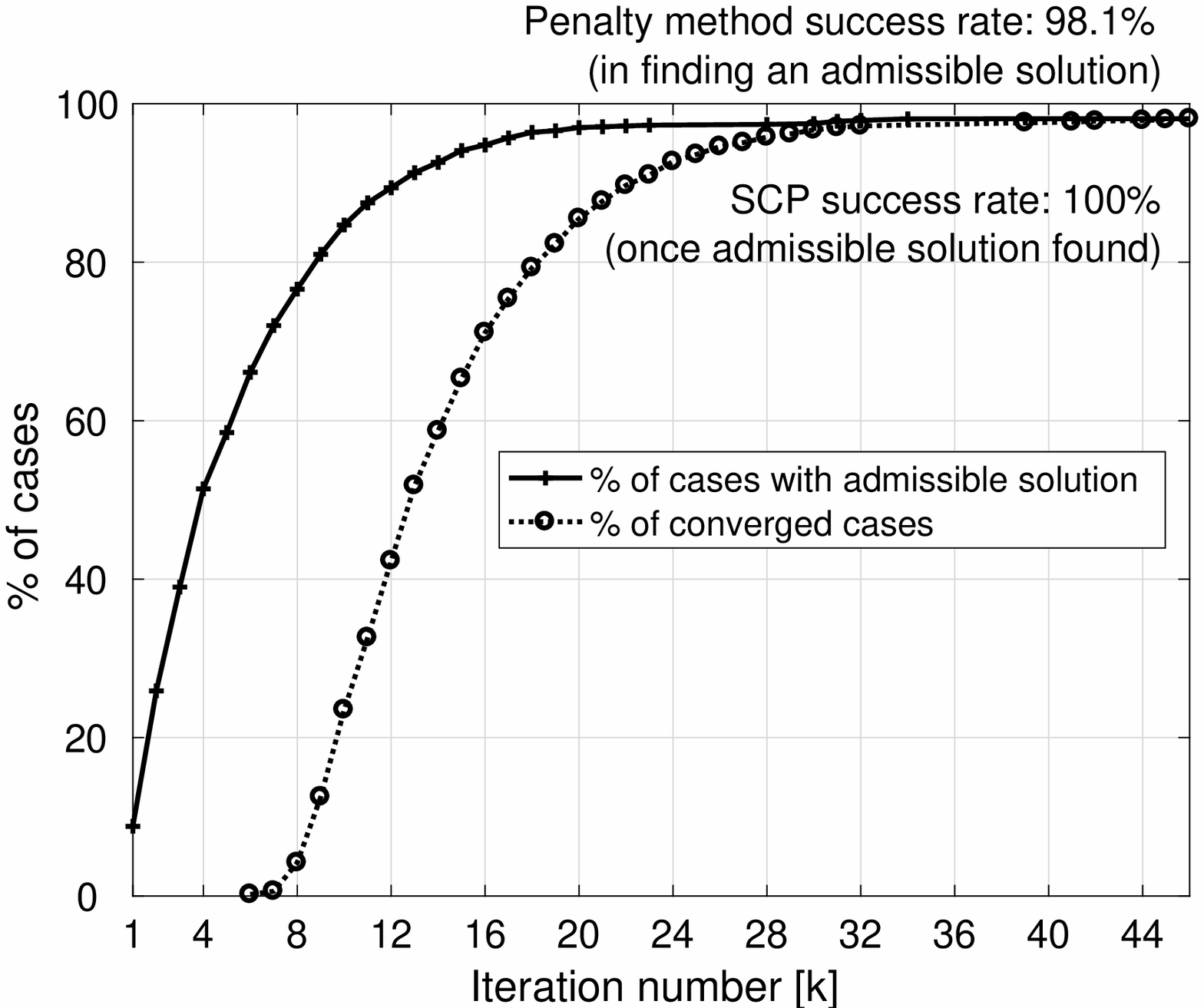}
	\caption{Convergence speed.\label{fig:converged_cdf}}
\end{figure}

Figure~\ref{fig:converged_cdf} shows the convergence of the proposed method on this illustrative problem. The penalty method is able to find an admissible solution in 98.1\% of the cases. Once an admissible solution is found the SCP enjoys guaranteed convergence, thus achieving a 100 \% success rate.

\begin{figure}
	\includegraphics[width=0.5\textwidth]{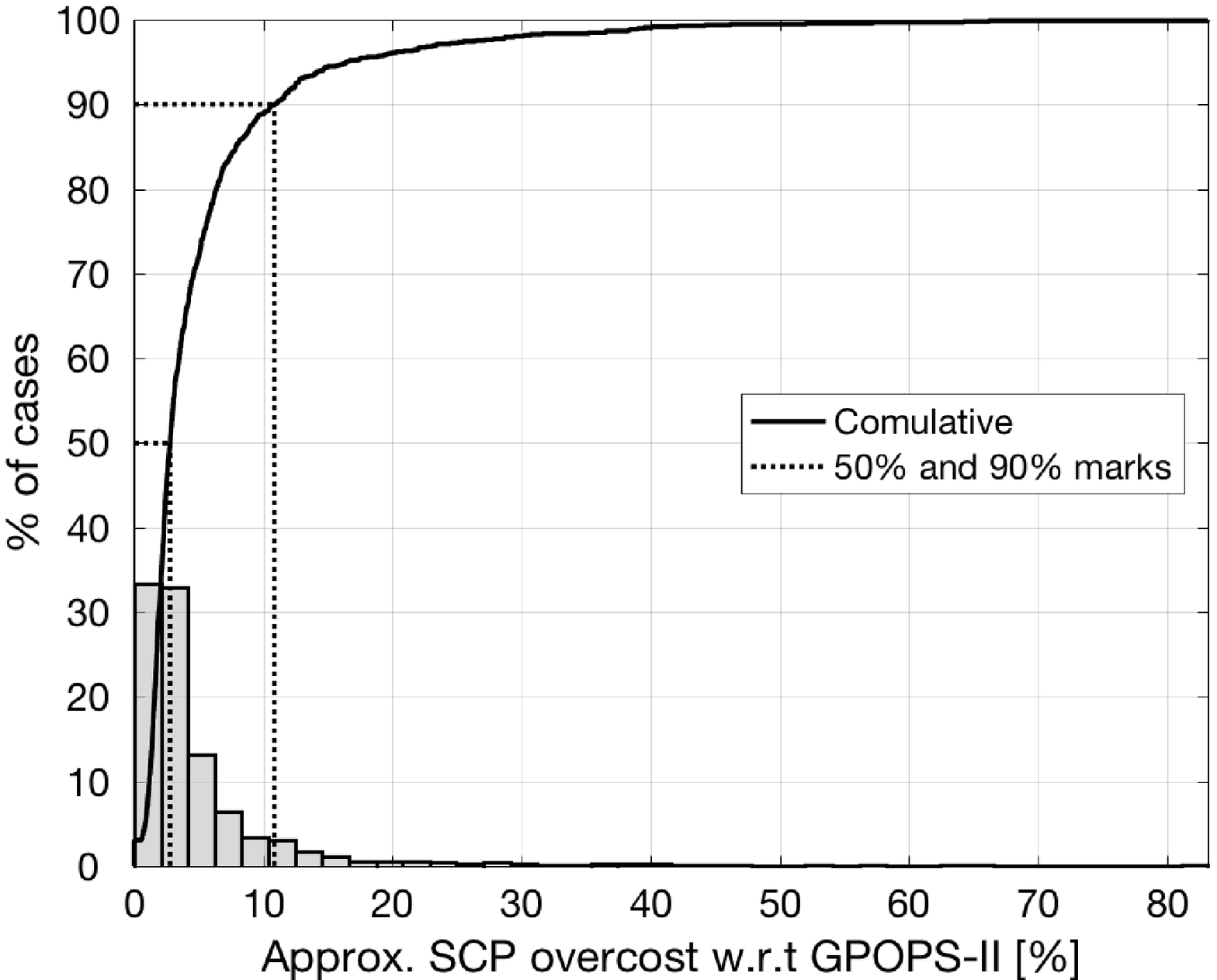}
	\caption{Comparison of the cost between the SCP and GPOPS-II.\label{fig:delta_cost_histogram}}
\end{figure}

Figure~\ref{fig:delta_cost_histogram} compares the cost of the solution offered by the proposed SCP with the solution obtained with the pseudospectral optimal control package GPOPS-II \cite{Patterson2014}. Here the GPOPS-II solution is used as a proxy for the globally optimal solution.  When the cost of both methods are compared it can be argued that the proposed SCP produces competitive solutions with only minor overcosts. For example, 50\% of the cases show an overcost below 5\%, while the overcost of 90\% of the cases is below 12\%.

Where the proposed method shows one of its most remarkable advantages is in its low computational cost. A single-threaded C implementation of the proposed SCP has been tested on a laptop and on two embedded computing platforms: an NVIDIA TX2 and Raspberry Pi 3 B+ (RPi). To solve the convex programming problems the interior-point solver IPOPT \cite{Wachter2005}, paired with HSL's MA27 linear solver routine \cite{HSL}, is used.

\begin{figure}
	\includegraphics[width=0.5\textwidth]{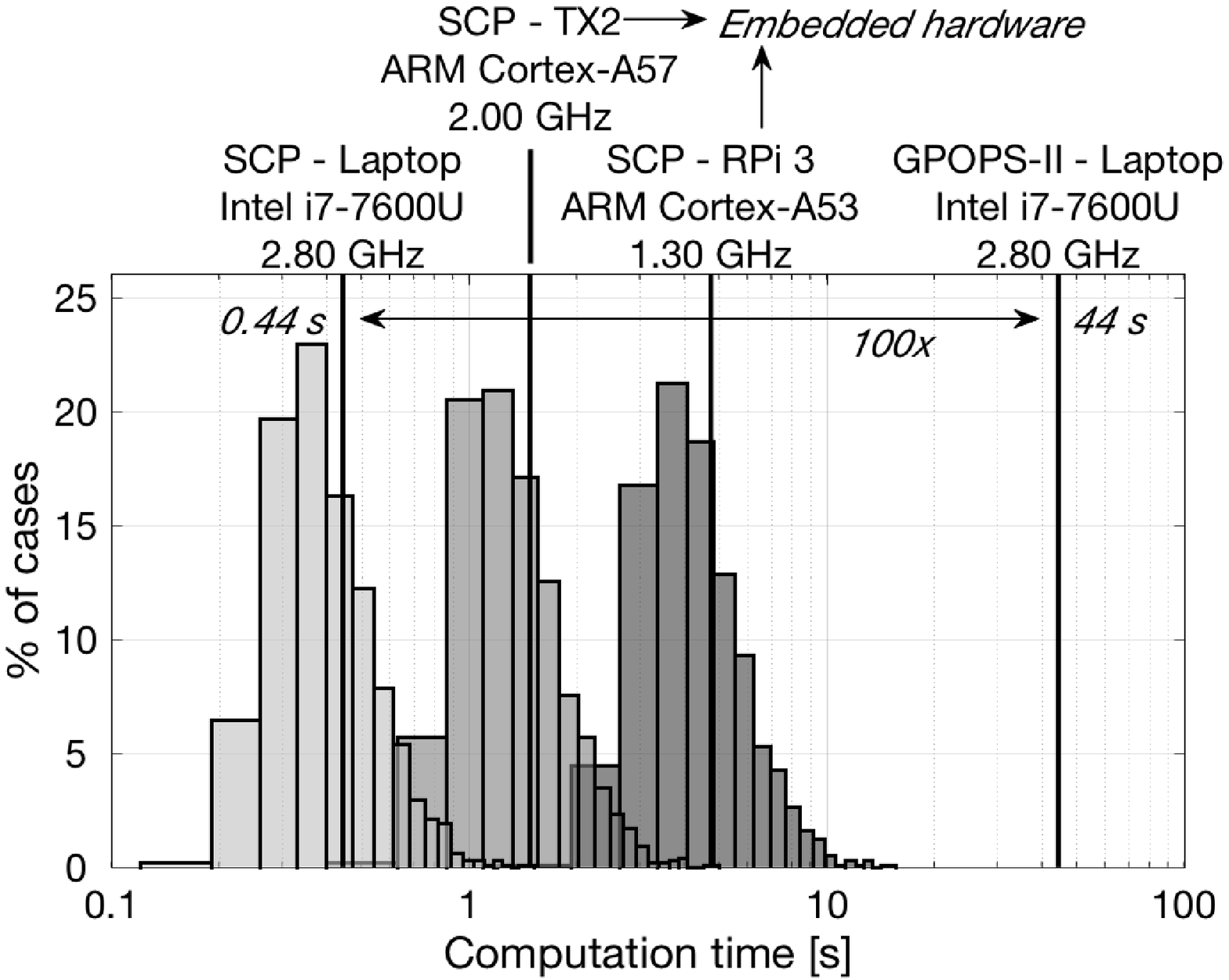}
	\caption{SCP's computational time on different computing platforms, compared to GPOPS-II.\label{fig:total_time_histogram}}
\end{figure}

Figure~\ref{fig:total_time_histogram} shows an histogram of the computational time required to solve the problems on these three platforms, with the solid lines indicating the averages. For reference, the average time taken by GPOPS-II on the laptop is also shown. The GPOPS-II solution uses the default settings (with IPOPT used as the underlying nonlinear programming solver). On a laptop machine the proposed SCP is shown to be a 100 times faster than GPOPS-II. Additionally the results on the two consumer-grade embedded computers provide empirical evidence that the proposed approach is suitable for implementation on embedded computers and fast enough to be used for real-time applications.  

\begin{figure}
	\includegraphics[width=0.5\textwidth]{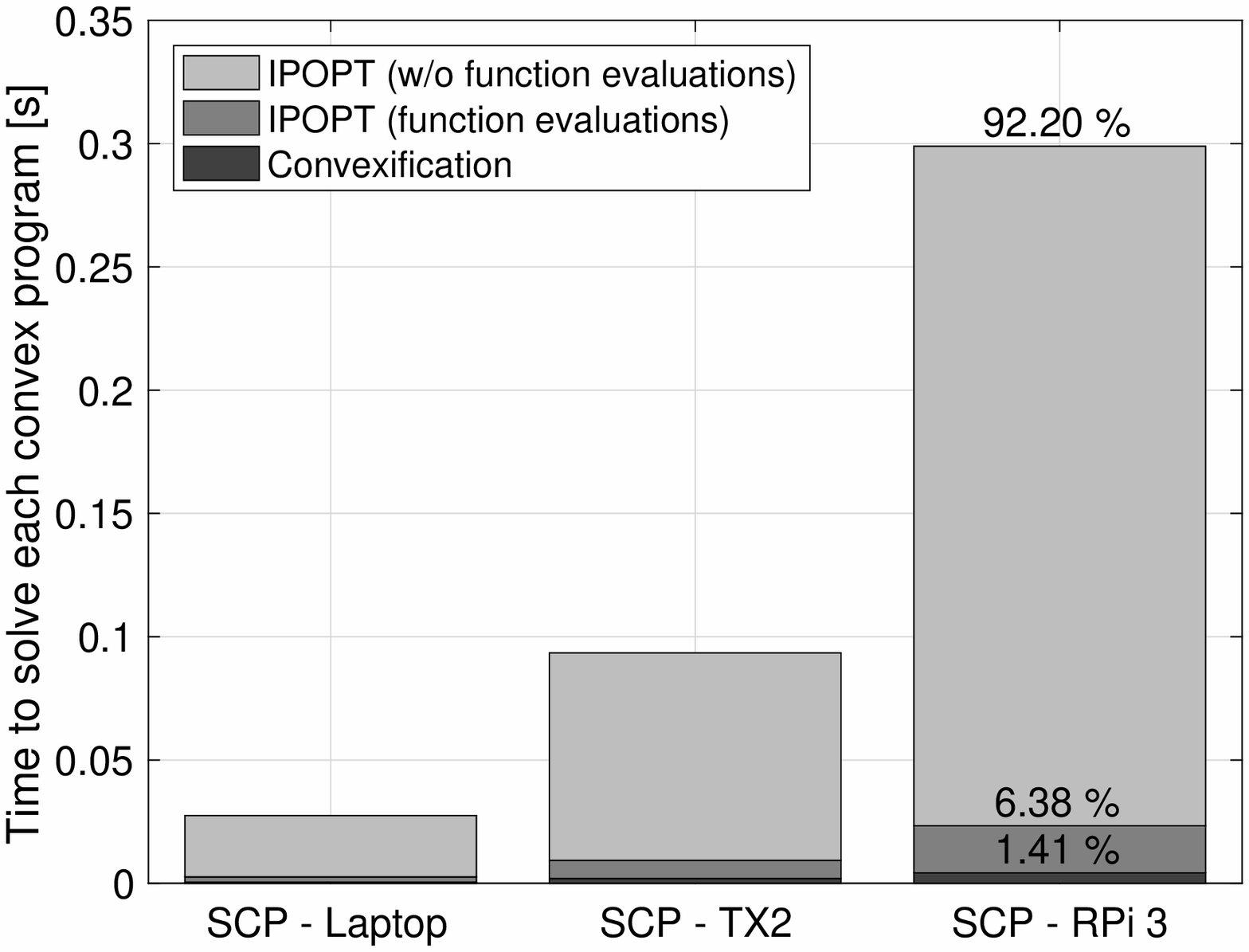}
	\caption{Average time to solve a single convex programming problem within the SCP sequence.\label{fig:iter_time}}
\end{figure}

Finally, Fig~\ref{fig:iter_time} shows how the time is spent solving each convex program within the SCP. Percentages for each of the major tasks are approximately the same across all platforms, but for readability they are only explicitly given for the RPi. What Fig~\ref{fig:iter_time} shows is that the proposed Taylor-based approach to obtain inner-convex approximations is computationally efficient, taking less than 1.5\% of the time per iteration. Also the time taken to evaluate the nonlinear functions (\emph{i.e.}, Eq.~\eqref{eq:gcvx}) is small, just below 6.5\%. The bulk of the time, $> 92\%$, is used to solve a linear system at each iteration of the IPOPT's interior-point method.

\subsection{Applicability to other problems}
The proposed SCP and Taylor-based convexification approach can be applied to a wide range of problems. The example provided here is only one of these potential applications. The authors have  applied this optimization approach to other problems, namely to the spacecraft attitude maneuver \cite{Virgili_2018_ASC_ConvexAttitude}.


\section{Concluding remarks}
When inner-convex approximations are used, a recursively feasible, descending, and convergent sequential convex programming (SCP) procedure is obtained for the wide class of non-convex programming problems with linear equality constraints.
This SCP procedure avoids the common pitfalls of infeasible iterates and cost increases of linearization-based SCP procedures without the need of trust regions.
By convexifiying each term of the Taylor series expansion of a non-convex function we obtain a computationally efficient and widely applicable method to generate inner-convex approximations. These Taylor-based approximations preserve the positive semidefinite part of the Hessian and thus provide a locally quadratic rate of convergence when the local minima is non-degenerate.
Given the attractive computational properties of the proposed approach it appears that the proposed SCP is suitable for onboard implementation and real-time use. The results obtained when solving a trajectory optimization problem on two embedded computers empirically substantiate this claim.


\bibliographystyle{elsarticle-num}
\bibliography{Biblio}

\end{document}